\documentclass{amsart}

\theoremstyle{definition}

\theoremstyle{remark}

\reversemarginpar

\begin{document}

\title{A dozen integrals: Russell-style}

\author{Tewodros Amdeberhan}
\address{Department of Mathematics,
MIT, Cambridge, MA 02139}
\email{tewodros@math.mit.edu}

\author{Victor H. Moll}
\address{Department of Mathematics,
Tulane University, New Orleans, LA 70118}
\email{vhm@math.tulane.edu}

%    General info
\subjclass{Primary 33}

\date{\today}


\maketitle

\newcommand{\ba}{\begin{eqnarray}}
\newcommand{\ea}{\end{eqnarray}}
\newcommand{\ift}{\int_{0}^{\infty}}
\newcommand{\ione}{\int_{0}^{1}}
\newcommand{\ifft}{\int_{- \infty}^{\infty}}
\newcommand{\no}{\noindent}
\newcommand{\realpart}{\mathop{\rm Re}\nolimits}
\newcommand{\imagpart}{\mathop{\rm Im}\nolimits}

\newtheorem{Definition}{\bf Definition}[section]
\newtheorem{Thm}[Definition]{\bf Theorem} 
\newtheorem{Example}[Definition]{\bf Example} 
\newtheorem{Lem}[Definition]{\bf Lemma} 
\newtheorem{Note}[Definition]{\bf Note} 
\newtheorem{Cor}[Definition]{\bf Corollary} 
\newtheorem{Prop}[Definition]{\bf Proposition} 
\newtheorem{Problem}[Definition]{\bf Problem} 
%\numberwithin{equation}{section}

On June 15, 1876, the Proceedings of the Royal Society of London 
published the paper \cite{russell1} by Mr. W. H. L. Russell entitled 
{\em On certain integrals}. The paper starts with {\em The following are 
certain integrals which will, I hope, be found interesting}. The rest of the 
paper is a list of $12$ definite integrals starting with
\begin{equation}
\ift dz \, e^{-(r+1)z+xe^{-z}} = e^{x}
\left( 1 -r/x + r(r-1)/x^{2} - \cdots \right)/x. 
\nonumber
\end{equation}

In honor of Russell notable achievements in the evaluation of integrals, we
hope the reader will find the next list interesting: \\

\begin{equation}
\ift x \left( \frac{\gamma \sinh \gamma x}{\cosh^{2} \gamma x} 
e^{-x^{2}/\pi^{2}} + \frac{\sqrt{\pi} \sinh x}{\cosh^{2}x} 
e^{- \gamma^{2} x^{2}} \right) \, dx = 
\ift \frac{e^{-x^{2}/\pi^{2}} \, dx}{\cosh \gamma x}  
\end{equation}
\begin{equation}
\ift x \left( \frac{1}{\pi} e^{-x^{2}/\pi^{2}} + \frac{1}{\sqrt{\pi}} 
e^{-x^{2}} \right) \frac{\sinh x \, dx}{\cosh^{2}x}  = 
\ift \frac{e^{-x^{2}} \, dx}{\cosh \pi x}  
\end{equation}
\begin{equation}
\ift x \left(  e^{-x^{2}/\pi} + 2  
e^{-16x^{2}/\pi} \right) \frac{\sinh 2 x \, dx}{\cosh^{2}2x}  = 
\ift \frac{e^{-4x^{2}/\pi} \, dx}{\cosh 4 x}  
\end{equation}
\begin{equation}
\int_{0}^{\infty} \left( \frac{\sinh x}{\cosh^{2}x} + 
\frac{\pi^{3/2} \sinh \pi x}{\cosh^{2} \pi x} \right) xe^{-x^{2}} \, dx  = 
\ift \frac{e^{-x^{2}} \, dx}{\cosh x}  
\end{equation}
\begin{equation}
\ift xe^{-x^{2}/\pi} \frac{\sinh x \, dx}{\cosh^{2}x}  =
\ift \frac{e^{-4x^{2}/\pi} \, dx}{\cosh 2x}  
\end{equation}
\begin{equation}
\int_{0}^{1} \frac{x^{-\ln x} \, dx}{1+x^{2}}   = 
\ift \frac{e^{-4x^{2}/\pi} \, dx}{\cosh 2 \sqrt{\pi}x } \, dx 
\end{equation}
\begin{equation}
\ift \frac{x^{2} e^{-x^{2}} \, dx}{\cosh \sqrt{\pi} x}  = 
\frac{1}{4} \ift \frac{e^{-x^{2}} \, dx}{\cosh \sqrt{\pi} x} 
\end{equation}
\begin{equation}
\ift \left( e^{-x^{2}/\pi^{2}} + \pi^{5/2} e^{-x^{2}} \right) 
\frac{x^{2} \, dx}{\cosh x}   =  
\frac{\pi^{2}}{2} \ift \frac{e^{-x^{2}/\pi^{2}} \, dx}{\cosh x}  
\end{equation}
\begin{equation}
\ift \left( \sqrt{\pi} e^{-x^{2}/3} + 9 \sqrt{3} \pi^{-2} e^{-3x^{2}/\pi^{2}} 
\right) \frac{x^{2} \, dx}
{\cosh x}  
 =  \frac{3 \pi \sqrt{3}}{2} 
\ift \frac{e^{-3x^{2}} \, dx}{\cosh \pi x} 
\end{equation}
\begin{equation}
\ift \left( \pi^{5} e^{-\pi^{3}x^{2}/G} + G^{5/2} e^{-Gx^{2}/\pi} \right)
\frac{x^{2} \, dx}{\cosh \pi x}  =  \frac{\pi G^{3/2}}{2}
\ift \frac{e^{-Gx^{2}/\pi} \, dx}{\cosh \pi x}    
\end{equation}
\begin{equation}
\ift x(3- 4 \pi x^{2} ) \frac{e^{-\pi x(x+1)} \, dx}{\sinh \pi x}  = 
\frac{1}{2 \pi} 
\end{equation}
\begin{equation}
\ift \frac{\sin^{2}x}{\cosh x + \cos x} \frac{dx}{x^{2}} + 
\frac{2}{\pi} \int_{0}^{e^{-\pi/2}} \frac{\tan^{-1}x}{x} \, dx = \frac{\pi}{4}
\end{equation}

\medskip

Here, 
\begin{equation}
\gamma := \lim\limits_{n \to \infty} \sum_{k=1}^{n} \frac{1}{k} - \ln n
\end{equation}
\noindent
is Euler's constant and 
\begin{equation}
G := \sum_{k=0}^{\infty} \frac{(-1)^{k}}{(2k+1)^{2}},
\end{equation}
\noindent
is Catalan's constant. We note the fact that both of these constants have 
so far resisted all attempts at proofs of irrationality. 

\medskip

Every formula can be checked by observing that if $f(y) = 1/\cosh(y)$  and 
$M$ is the transformation 
\begin{equation}
M(f)(y) = \ift e^{-x^{2}} f(xy) \, dx,
\nonumber
\end{equation}
\noindent
then we have the elementary relation
\begin{equation}
\frac{df}{dx} = - x f(y) \sqrt{1 - f^{2}(y)} 
\nonumber
\end{equation}
\noindent
and also 
\begin{equation}
y M(f)(y) = \sqrt{\pi} M \left( f \left( \frac{\pi}{y} \right) \right). 
\nonumber
\end{equation}
\noindent
Applying these to the left hand side integrands  produces the right-hand side.
In the last few calculations, the role of the 
functions $1/\sinh y$ and $e^{-x^{2}}$  has to be interchanged. \\

\noindent
{\bf Acknowledgements.} The second author acknowledges the partial support of 
NSF-DMS 04099658.

\end{document}